\documentclass[12]{article}
\usepackage{amssymb}
\usepackage{latexsym}
\usepackage{amsmath}
\usepackage{mathrsfs}
\usepackage[all]{xy}
\usepackage{enumerate}
\usepackage{amsthm}
\usepackage[dvips]{graphicx}
\usepackage{psfrag}
\usepackage{ifthen}
\usepackage[usenames]{color}

\newtheorem{theorem}{Theorem}[section]
\newtheorem{lemma}[theorem]{Lemma}
\newtheorem{metalemma}[theorem]{Meta Lemma}
\newtheorem{proposition}[theorem]{Proposition}
\newtheorem{corollary}[theorem]{Corollary}

\theoremstyle{definition}
\newtheorem{definition}[theorem]{Definition}

\newtheorem{example}[theorem]{Example}

\newcommand{\comment}[1]
	   {\ifthenelse{\equal{\showcomments}{yes}}
	     {\footnotemark\marginpar{\sffamily{\tiny
		   \addtocounter{footnote}{-1}\footnotemark#1
}\normalfont}}{}}

\newcommand{\tr}{\textrm}
\newcommand{\showcomments}{yes}

\newcommand{\mo}{{-1}}
\newcommand{\bel}[1]{\begin{equation}\label{#1}}
\newcommand{\ee}{\end{equation}}

\newcommand{\bdy}{\partial}

\newcommand{\ol}{\overline}

\newcommand{\LBA}{\left\{\begin{array}}
\newcommand{\EAR}{\end{array}\right.}
\input epsf

\title{Linear estimates for solutions of quadratic equations in free groups}
\author{Olga Kharlampovich and Alina Vdovina}

\begin{document}
\maketitle

\begin{abstract}
We prove that in a free group
the length of the value of each variable in a minimal  solution of  a standard  quadratic equation is bounded by $2s$ for orientable equation and by $12s^4$ for non-orientable equation, where $s$ is the sum of the lengths of the coefficients.
\end{abstract}

\section{Introduction}

The study of quadratic equations over free groups started with the
work of Malc$'$ev \cite{Malcev-1962} and has been deepened extensively
ever since.  One of the reasons research in this topic has been so
fruitful is the deep connection between quadratic equations and the
topology of surfaces.

In \cite{Comerford-Edmunds-1981} the problem of deciding if a quadratic
equation over a free group is satisfiable
was shown to be decidable. In
addition it was shown in \cite{Olshanskii-1989},
\cite{Grigorchuk-Kurchanov-1989}, and \cite{Grigorchuk-Lysenok-1992}
that if $n$, the number of variables, is fixed, then deciding if a standard
quadratic equation has a solution can be done in time which is polynomial in
the sum of the lengths of the coefficients.  These results imply that
the problem is solvable in at most exponential time. In
\cite{KLMT} it was shown that the problem of deciding if a quadratic
equation over a free group is satisfiable is NP-complete.
 We will improve
on this by proving  that  in a free group,
the length of the value of each variable in a minimal solution of a standard  quadratic equation is bounded by $2s$ for orientable equation and by $3s^2$ for non-orientable equation, where $s$ is the sum of the lengths of the coefficients.

\medskip

\par {\bf Definition 0.} Let $G$ be a group and $w$ be an element of its
commutator subgroup. We define the {\em orientable genus\/} $g(w)$ of
 $w$ as the least positive integer $g$ such that $w$ is a product
 of $g$ commutators in $G$.
 We define the {\em non-orientable genus\/} $g(w)$ of
 $w$ as the least positive integer $n$ such that $w$ is a product
 of $g$ squares in $G$.
\medskip

\par {\bf Theorem 1.} Let $C$ be an orientable (resp., non-orientable) word of
 genus $g$ in a free group $F$. Then $C$
can be presented in the form $[a_1,b_1]\ldots[a_g,b_g]$ (resp. $d_1^2\ldots d_k^2[a_{1},b_{1}]\ldots[a_m,b_m]$ with $k+2m=g$),
where $|a_i|< 2|C|$,$|b_i| < 2|C|$,$|d_i|<2|C|$ for $i=1,...,g$.

If $C$ is non-orientable then $C$ can be represented
as a product of squares $a_1^2...a_g^2$ with $|a_i| \leq 3|C|^2$.
\medskip

 An {\it orientable quadratic set of words}  is a set of cyclic  words
$ w_1,w_2,\dots ,w_{k}$
 (a cyclic word is the orbit of a linear word under cyclic permutations)
 in some alphabet $a_1^{\pm 1},a_2^{\pm 1},\dots$ of letters
 $a_1,a_2,\dots$ and their inverses $a_1^{-1},a_2^{-1},\dots$ such that
\begin{itemize}
\item[(i)] if $a_i^\epsilon$ appears in $w_i$ (for $\epsilon\in\{\pm 1\}$)
 then $a_i^{-\epsilon}$ appears exactly once in $w_j$,
\item[(ii)] the word $w_i$ contains no cyclic factor (subword of
 cyclically consecutive letters in $w$) of the form $a_l a_l^{-1}$ or
 $a_l^{-1}a_l$ (no cancellation),

\end{itemize}
To define a {\it quadratic set of words} we replace condition (i) by the condition that each $a_i^{\pm 1}$ appears in the set of words $ w_1,w_2,\dots ,w_{k}$ exactly twice.
Quadratic sets of words were first defined in \cite{LS}, p.60.

The {\it genus}  of a quadratic set of words is defined as the sum of genera of the surfaces obtained from
 $k$ discs with words $ w_1,w_2,\dots ,w_{k}$ on their boundaries when we identify the edges labeled by the same letters.
\medskip

The Proposition below follows from Olshanskii's theorem (Theorem 4), described later:

\par{\bf Proposition } The following two conditions on a set $\{C_1,...,C_k\}$
of elements of a free group $F$ are equivalent:

(a) {\it The system $\{W_i=C_i, i=1,...,k\}$ has a solution in $F$ where $\{W_i\}$ (orientable or non-orientable) quadratic set of words of genus $g$};

and

(b){\it The standard (resp. orientable or non-orientable quadratic equation of genus $g$ with coefficients  $\{C_1,...,C_k\}$ has a solution in $F$.}

{\bf Proof} By Olshanskii's theorem, (b) implies that there is a collection of discs $D_1,\ldots , D_m$ with boundaries labelled by a quadratic set of words $W_1,\ldots ,W_m$ in some variables  $P=\{p_1,\ldots ,p_n\},$ and
there is a mapping  $\ol{\psi}:P \rightarrow (A\cup A^\mo)^*$ such
    that upon substitution, the coefficients $C_1,\ldots,C_{m-1}$ and
    $C$ can be read without cancellations around the boundaries of
    $D_1,\ldots,D_{m-1}$ and $D_m$, respectively. Then $\ol{\psi}(p_1),\ldots ,\ol{\psi}(p_n)$ is a solution of the system in $W_1=C_1,\ldots ,W_m=C.$
    Let $g_0,\ldots, g_l$ be genera of the surfaces obtained from
 $m$ discs with words $ W_1,W_2,\dots ,W_{k}$ on their boundaries when we identify the edges labeled by the same letters.
Inequalities in (iii) imply that $\sum _{i=0}^lg_i\leq g.$ If this inequality is strict, and $g-\sum _{i=0}^lg_i=r>0$, we can consider instead of $W_1=C_1$ equation  $W_1[p_{n+1},p_{n+2}]\ldots [p_{n+2r-1},p_{n+2r}]=C_1$  and define $\ol{\psi}(p_j)=1$ for $j=n+1,\ldots ,n+2r.$ Then
$W_1[p_{n+1},p_{n+2}]\ldots [p_{n+2r-1},p_{n+2r}], W_2,\ldots ,W_m$ has genus $g$. This proves (a).

Now, suppose we have (a), and $W_i$ are equations in variables $P$.  If some letter $p_i\in P$ is contained in different equations we can express it from one equation and substitute in the other.
For example, $W_1=W_{11}p_1W_{12}=C_1$ and $W_2=W_{21}p_1^{-1}W_{22}=C_2$ become one equation $$W_{21}W_{12}C_1^{-1}W_{11}W_{22}=C_2$$ which can be rewritten as
$ W_{21}W_{12}W_{11}W_{22}=C_1^{(W_{21}W_{12})^{-1}}C_2.$ If we label the edges of two polygons $D_1$ and $D_2$ by $W_1$ and $W_2$, then the left hand side of this equation will be written on the boundary of a  polygon obtained  by identifying the edges labeled by $p_1$ and removing them. We continue this procedure until there is no letter $p_i$ contained in two different equations. Taking the inverses of both sides of each  equation, we obtain a system $L_i(P)=R_i, j=1,\ldots ,s$ such that both appearances of every letter from $P$ are contained in one word $L_i$ and each $R_i$ is a product of some conjugates of coefficients $C_j^{-1}$ (if the set of words is orientable) or $C_j^{\pm 1}$ (if the set of words is non-orientable).
The sum of topological genera (see definition in Section 3) of words $L_i$ is the
sum of genera of the surfaces obtained from $s$ discs with words $L_1,\ldots ,L_s$ on their boundaries when we identify the edges labeled by the same letters. By construction, the same surfaces are obtained  from $k$ discs with words $W_1,\ldots ,W_k$ on their boundaries. The sum of their genera is $g$. Topological and algebraic genus of each word $L_i$ is the same. Therefore, the left side of the equation
$$L_1\ldots L_s= R_1\ldots R_s$$ has algebraic genus $g$. In the orientable case this proves (b). In the non-orientable case we have to change variables and replace negative powers of $C_i$'s by positive powers. Changing variables we can assume that the left-hand side is a product of $g$ squares. We can also assume that the right-hand side is a product of two parts. The first part is the product of conjugates of negative powers of $C_i$'s and  the second part is the product of conjugates of positive powers of $C_i$'s. Therefore $$L_1\ldots L_s= R_1\ldots R_s$$ can be transformed to the form $x_1^2\ldots x_g^2=P_1P_2,$ where $P_1=\prod C_{i_j}^{-Z_j}$ and $P_2=\prod C_{i_l}^{Z_l}$. We will write it as
$$x_1^2\ldots x_gP_2^{-1}P_2x_gP_2^{-1}=P_1,$$ make a substitution $\bar x_g=x_gP_2^{-1}$ and re-write as
$$x_1^2\ldots \bar x_g^2P_2^{\bar x_g}P_1^{-1}=1.$$ Now conjugates of all $C_i$'s appear  in positive exponents, conjugating again we can put them in the right order.  Since the system in (a) has a solution, the standard quadratic equation also has a solution,  and (b) is proved.

\medskip
\par {\bf Theorem 2.} Let $W_1$,...,$W_k$ be an orientable quadratic set of words
of  genus $g$, and $C_1,\ldots, C_k$ be elements of a free group $F$ such that the system $W_i=C_i,\ i=1,\ldots, k$ has a solution in $F$  and $\sum_{i=1}^k|C_i|=s$ in $F$.
 Then some product of conjugates of $C_i$'s in any order
can be presented as a product of at most $g$ commutators  of elements in $F$ with lengths strictly
less then $2s$, and conjugating elements also have length bounded by $2s$.

Let $W_1$,...,$W_k$ be a non-orientable quadratic set of words
of  genus $g$, and $C_1,\ldots, C_k$ be elements of a free group $F$ such that the system $W_i=C_i,\ i=1,\ldots, k$ has a solution in $F$  and $\sum_{i=1}^k|C_i|=s$ in $F$.
Then some product of conjugates of $C_i$'s in any order
can be presented as a product of at most $g$ squares $a_1^2...a_g^2$ with $|a_i| \leq 12s^4$
and conjugating elements  have length bounded by $2s^2$.

\medskip
\par {\bf Theorem 3.} Let $h$ be an orientable word of
 genus $g$ in a hyperbolic group $\Gamma$.
Let $M$ be the number of elements in $\Gamma$ represented
by words of length at most $4\delta$ in $F(X)$ ($\delta$ is the hyperbolicity constant), $l=\delta(log_2(12g-6)+1)$.
Then $h$
can be presented in a form $[a_1,b_1]\ldots[a_g,b_g]$,
where $|a_i|< 2|h|+3(12g-6)(12l+M+4)$, $|b_i|< 2|h|+3(12g-6)(12l+M+4)$ for $i=1,...,g$.

\medskip

\section{Quadratic equations}
A quadratic equation $E$ with variables $\{x_i,y_i,z_j\}$ and
non-trivial coefficients $\{C_i,C\} \in F(A)$ is said to be in
\emph{standard form} if its coefficients are expressed as freely and
cyclically reduced words in $A^*$ and $E$ has either the form:
\bel{eqn:orientable} \left( \prod_{i=1}^g [x_i,y_i]
\right)\left(\prod_{j=1}^{m-1} z_j^\mo C_j z_j\right) C = 1 \tr{~or~}
\left(\prod_{i=1}^g [x_i,y_i]\right)C=1\ee where $[x,y]=x^\mo y^\mo
xy$, in which case we say it is \emph{orientable} or it has the form
\bel{eqn:non-orientable} \left(\prod_{i=1}^g x_i^2\right)
\left(\prod_{j=1}^{m-1} z_j^\mo C_j z_j\right) C = 1 \tr{~or~}
\left(\prod_{i=1}^g x_i^2\right)C=1\ee in which case we say it is
non-orientable.  The \emph{genus} of a quadratic equation is the
number $g$ in (\ref{eqn:orientable}) and (\ref{eqn:non-orientable})
and $m$ is the number of coefficients. If $g=0$ then we will define
$E$ to be orientable. If $E$ is a quadratic equation we define its
\emph{reduced Euler characteristic}, $\ol{\chi}$ as follows:
\[\ol{\chi}(E) =\left\{
  \begin{array}{l}2-2g \tr{~if~} E \tr{~is orientable} \\ 2-g
    \tr{~if~} E \tr{~is not orientable}\end{array}\right.\] We finally
  define the \emph{length} of a quadratic equation $E$ to be\[
  \tr{length}(E) = |C_1|+\ldots+|C_{n-1}|+C+ 2\tr{(number of
    variables)}\] It is a well known fact that an arbitrary quadratic
  equation over a free group can be brought to a standard form in time
  polynomial in its length.

\par {\bf Corollary 1 of Theorem 2.} Let $E(X)=1$ be a standard consistent quadratic equation in a free group $F$ with the set of variables $X$.
For the orientable equation there exists a solution $\phi :F(X)*F\rightarrow F$ such that for each $x\in X$,
$|\phi (x)|\leq 2s,$ where $s$ is the sum of the lengths of the coefficients.
For the non-orientable equation, $|\phi (x)|\leq 12s^4,$

Let $G$ be a group. Let $G[X] = G \ast F(X)$, where $X = \{x_1,\ldots, x_n\}$.
Let $S(X) = 1$ be a system of equations over $G$, that is, $S \subset G[X]$.
By $V_G(S)$ denote the set of all solutions in $G$ of the system $S(X) = 1$, it is called the
algebraic set defined by $S$.
$V_G(S)$ uniquely corresponds to the normal subgroup
$$R(S) = \{T(x) \in G[X] \mid \forall g \in G^n\ (S(g) = 1 \rightarrow T(g) = 1) \}$$
of the group $G[X]$.
The quotient group
$$G_{R(S)} = G[X]/R(S)$$ is called the coordinate group of the system $S(X)=1$.

\vspace{2mm}
Let $G$ be a group with a generating set $A$. A system of equations
$S = 1$  is called {\em triangular quasi-quadratic} (shortly, TQ)
over $G$ if it can be partitioned into the following subsystems

\medskip
$S_1(X_1, X_2, \ldots, X_n,A) = 1,$

\medskip
$\ \ \ \ \ S_2(X_2, \ldots, X_n,A) = 1,$

$\ \ \ \ \ \ \ \ \ \  \ldots$

\medskip
$\ \ \ \ \ \ \ \ \ \ \ \ \ \ \ \ S_n(X_n,A) = 1$

\medskip \noindent
 where for each
$i$ one of the following holds:
\begin{enumerate}
\item [1)] $S_i$ is quadratic  in variables $X_i$;
 \item [2)] If we denote  $G_{i}=G_{R(S_{i}, \ldots, S_n)}$ for
$i = 1, \ldots, n$, and put $G_{n+1}=G.$ $S_i= \{[y,z]=1, [y,v]=1, v\in C(u) \mid y, z \in X_i\}$ where $u$ is a
group word in $X_{i+1} \cup  \ldots \cup X_n \cup A$, and $C(u)$ the centraliser of $u$ in $G_{i+1}$. In this case
we say that $S_i=1$ corresponds to an extension of a centraliser;
 \item [3)] $S_i= \{[y,z]=1 \mid y, z \in X_i\}$;
 \item [4)] $S_i$ is the empty equation.
  \end{enumerate}

The  TQ system $S = 1$ is
called {\em non-degenerate} (shortly, NTQ) if the following
condition hold:
 \begin{enumerate}
  \item [5)]  each system $S_i=1$, where $X_{i+1}, \ldots, X_n$
are viewed as the corresponding constants from $G_{i+1}$ (under the
canonical maps $X_j \rightarrow G_{i+1}$, $j = i+1, \ldots, n$)  has
a solution in $G_{i+1}$;

  \end{enumerate}

\par {\bf Corollary 2 of Theorem 2.}  Let $E(X)=1$ be an NTQ system in the standard form (all quadratic equations are in the standard form) in a free group $F$ with the set of variables $X$ and  with $l$ levels. Let $s_i$ be the total length of the coefficients on the level $i$.
There exists a solution $\phi :F(X)*{F}\rightarrow {F}$ such that for each $x\in X$,
$|\phi (x)|\leq 12^{2l}s_1^4\ldots s_l^4.$

\section{Orientable and non-orientable Wicks forms}

 {\bf Definition 1.}
 An {\it orientable Wicks form\/} is a cyclic word $w$
 (a cyclic word is the orbit of a linear word under cyclic permutations)
 in some alphabet $a_1^{\pm 1},a_2^{\pm 1},\dots$ of letters
 $a_1,a_2,\dots$ and their inverses $a_1^{-1},a_2^{-1},\dots$ such that
\begin{itemize}
\item[(i)] if $a_i^\epsilon$ appears in $w$ (for $\epsilon\in\{\pm 1\}$)
 then $a_i^{-\epsilon}$ appears exactly once in $w$,
\item[(ii)] the word $w$ contains no cyclic factor (subword of
 cyclically consecutive letters in $w$) of the form $a_i a_i^{-1}$ or
 $a_i^{-1}a_i$ (no cancellation),
\item[(iii)] if $a_i^\epsilon a_j^\delta$ is a cyclic factor of $w$ then
 $a_j^{-\delta}a_i^{-\epsilon}$ is not a cyclic factor of $w$
(substitutions of the form
 $a_i^\epsilon a_j^\delta\longmapsto x,
 \quad a_j^{-\delta}a_i^{-\epsilon}\longmapsto x^{-1}$ are impossible).
\end{itemize}

 An orientable Wicks form $w$ is an element of the commutator subgroup
 when considered as an element in the free group $G$ generated by
 $a_1,a_2,\dots$. We define the {\em algebraic genus\/} $g_a(w)$ of
 $w$ as the least positive integer $g_a$ such that $w$ is a product
 of $g_a$ commutators in $G$.

 The {\em topological genus\/} $g_t(w)$ of an orientable Wicks
 form $w$ is defined as the topological
 genus of the orientable compact connected surface obtained by
 labelling and orienting the edges of a $2e-$gon (which we
 consider as a subset of the oriented plane) according to
 $w$ and by identifying the edges in the obvious way.

{\bf Definition 2.}
 A {\it non-orientable Wicks form\/} is a cyclic word $w$
 (a cyclic word is the orbit of a linear word under cyclic permutations)
 in some alphabet $a_1^{\pm 1},a_2^{\pm 1},\dots$ of letters
 $a_1,a_2,\dots$ and their inverses $a_1^{-1},a_2^{-1},\dots$ such that
\begin{itemize}
\item[(i)] if $a_i^\epsilon$ appears in $w$ (for $\epsilon\in\{\pm 1\}$)
 then $a_i^{\pm \epsilon}$ appears exactly once in $w$, and there is at least one letter which appears with the same exponent,
\item[(ii)] the word $w$ contains no cyclic factor (subword of
 cyclically consecutive letters in $w$) of the form $a_i a_i^{-1}$ or
 $a_i^{-1}a_i$ (no cancellation),
\item[(iii)] if $a_i^\epsilon a_j^\delta$ is a cyclic factor of $w$ then
 $a_j^{-\delta}a_i^{-\epsilon}$ or another appearance of $a_i^\epsilon a_j^\delta$ is not a cyclic factor of $w$.

\end{itemize}

 A non-orientable Wicks form $w$ is an element of the subgroup of squares
 when considered as an element in the free group $G$ generated by
 $a_1,a_2,\dots$. We define the non-orientable {\em algebraic genus\/} $g_a(w)$ of
 $w$ as the least positive integer $g_a$ such that $w$ is a product
 of $g_a$ squares in $G$.

 The {\em topological genus\/} $g_t(w)$ of a non-orientable Wicks
 form $w$ is defined as the topological
 genus of the non-orientable compact connected surface obtained by
 labelling and orienting the edges of a $2e-$gon  according to
 $w$ and by identifying the edges in the obvious way.

{\bf Remark}
{\sl The algebraic and the topological genus of an orientable Wicks
 form coincide (cf. \cite{C,CE}).} The same is true for non-orientable Wicks forms.

 We define the {\em genus\/} $g(w)$ of a
 Wicks form $w$ by $g(w)=g_a(w)=g_t(w)$.

 Consider the orientable compact surface $S$ associated to an orientable (non-orientable)
 Wicks form $w$. This surface carries an immersed graph
 $\Gamma\subset S$ such that $S\setminus \Gamma$ is an open polygon
 with $2e$ sides (and hence connected and simply connected).
 Moreover, conditions (ii) and (iii) on Wicks form imply that $\Gamma$
 contains no vertices of degree $1$ or $2$ (or equivalently that the
 dual graph of $\Gamma\subset S$ contains no faces which are $1-$gones
 or $2-$gones). This construction works also
 in the opposite direction: Given a graph $\Gamma\subset S$
 with $e$ edges on an orientable (non-orientable) compact connected surface $S$ of genus $g$
 such that $S\setminus \Gamma$ is connected and simply connected, we get
 an orientable(non-orientable) Wicks form of genus $g$ and length $2e$ by labelling and
 orienting the edges of $\Gamma$ and by cutting $S$ open along the graph
 $\Gamma$. The associated orientable(non-orientable) Wicks form is defined as the word
 which appears in this way on the boundary of the resulting polygon
 with $2e$ sides. We identify henceforth orientable(non-orientable) Wicks
 forms with the associated immersed graphs $\Gamma\subset S$,
 speaking of vertices and edges of orientable (non-orientable) Wicks form.

 The formula for the Euler characteristic
 $$\chi(S)=2-2g=v-e+1$$ in orientable and $$\chi(S)=2-g=v-e+1$$ in non-orientable case
 (where $v$ denotes the number of vertices and $e$ the number
 of edges in $\Gamma\subset S$) shows that
 an orientable Wicks(non-orientable) form of genus $g$ has at least length $4g$($2g$)
 (the associated graph has then a unique vertex of degree $4g$
 and $2g$ edges) and at most length $6(2g-1)$($6(g-1)$) (the associated
 graph has then $2(2g-1)$($2(g-1)$) vertices of degree three and
 $3(2g-1)$($(3(g-1)$) edges).

 We call an orientable Wicks form of genus $g$ {\em maximal\/} if it has
 length $6(2g-1)$ in orientable and $6(g-1)$ in non-orientable case.

 A vertex $V$ (with oriented edges $a,b,c$ pointing toward $V$) is
 {\em positive\/} if
$$w=ab^{-1}\dots bc^{-1}\dots ca^{-1}\dots \quad {\rm or }\quad
 w=ac^{-1}\dots cb^{-1}\dots ba^{-1}\dots $$
 and $V$ is {\em negative\/} if
 $$w=ab^{-1}\dots ca^{-1}\dots bc^{-1}\dots \quad {\rm or }\quad
 w=ac^{-1}\dots ba^{-1}\dots ab^{-1}\dots \quad
 ..$$

\par Let $V$ be a negative vertex of an orientable maximal Wicks
form of genus $g>1$. There are three possibilities, denoted
configurations of type $\alpha,\ \beta$ and $\gamma$
(see Figure 1) for the local configuration at $V$.

\centerline{\epsfysize2.5cm\epsfbox{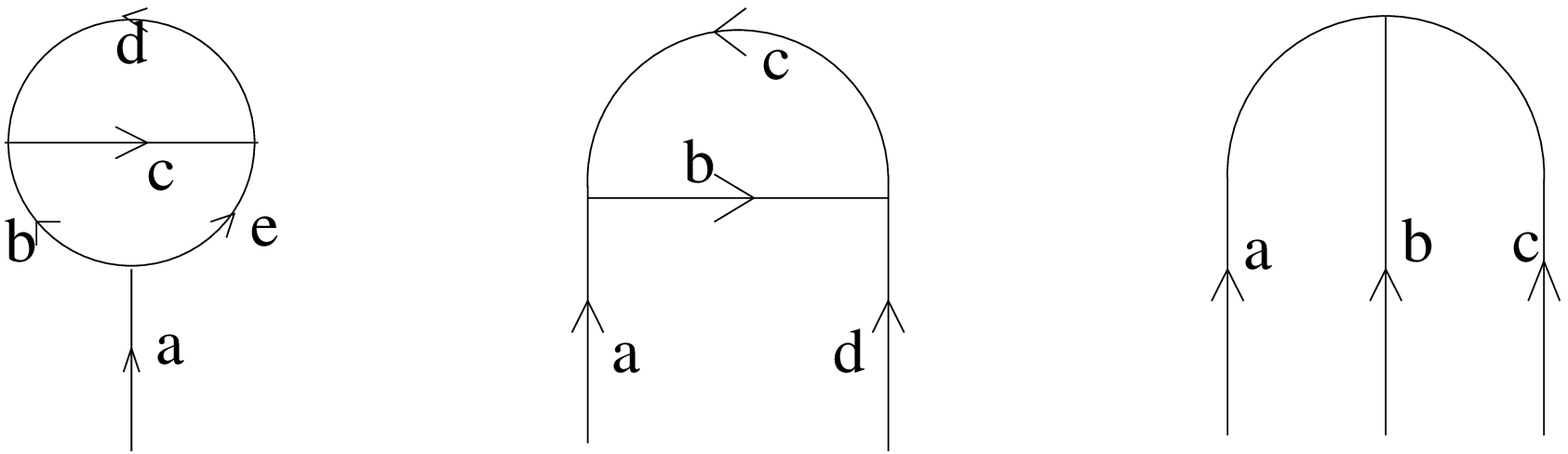}}
\centerline{\it Figure 1.}

\par Type $\alpha$. The vertex $V$ has only two neighbours
which are adjacent to each other. This implies that $w$ is of the
form
$$w=w_1x_1abcdb^{-1}ec^{-1}d^{-1}e^{-1}a^{-1}x_2w_2x_2^{-1}x_1^{-1}w_3$$
(where $w_1,w_2,w_3$ are subfactors of $w$)
and $w$ is obtained from the maximal orientable Wicks form
$$w'=w_1xw_2x^{-1}w_3$$
of genus $g-1$ by the substitution  $x\longmapsto
x_1abcdb^{-1}ec^{-1}d^{-1}e^{-1}a^{-1}x_2$ and $x^{-1}\longmapsto x_2^{-1}
x_1^{-1}$ (this construction is called the $\alpha -$construction in
{\cite{V}).

\par Type $\beta$. The vertex $V$ has two non-adjacent neighbours.
The word $w$ is then of the form
$$w=w_1x_1abca^{-1}x_2w_2y_1db^{-1}c^{-1}d^{-1}y_2w_3$$
(where perhaps $x_2=y_1$ or $x_1=y_2$, see \cite{V} for all the details).
The word $w$ is then obtained by a $\beta -$construction from the
word $w'=w_1xw_2yw_3$ which is an orientable maximal Wicks form
of genus $g-1$.

\par Type $\gamma$. The vertex $V$ has three distinct neighbours. We have
then
$$w=w_1x_1ab^{-1}y_2w_2z_1ca^{-1}x_2w_3y_1bc^{-1}z_2w_4$$
(some identifications among $x_i,\ y_j$ and $z_k$ may occur,
see \cite{V} for all the details) and the word $w$ is obtained by a
so-called $\gamma-$construction from the word $w'=w_1x w_3y
w_2 z w_4$, where $x=x_1x_2$, $y=y_1y_2$,$z=z_1z_2$.

For each of $\alpha$,$\beta$ or $\gamma$ an inverse transformation
 is well defined
(see \cite{V} for the details), so we can associate a Wicks form of genus $g-1$
to a Wicks form of genus $g$ and a given negative vertex.

\par {\bf Definition 3.}
We call the application which associates to an orientable maximal Wicks
form $w$ of genus $g$ with a chosen negative vertex $V$ the orientable
maximal Wicks form $w'$ of genus $g-1$ defined as above the {\it reduction}
of $w$ with respect to the negative vertex $V$.

\par An inspection of figure $1$ shows that reductions with respect to
vertices of type $\alpha$ or $\beta$ are always paired since two doubly
adjacent vertices are negative, of the same type
($\alpha$ or $\beta$) and yield the same reductions.

\par The above constructions of type $\alpha,\ \beta$ and $\gamma$
can be used for a recursive construction of all
orientable maximal Wicks forms of genus $g>1$ (\cite{V}).

Let $a$ be an edge of a non-orientable maximal genus $g$ Wicks form $w$ appearing with the same exponent.
Without loss of generality we can assume that $w=v_1av_2av_3$. The word $w^{\prime}$ obtained from $v_1v_2^{-1}v_3$
by simplifications and reductions is either non-orientable genus $g-1$ Wicks form or orientable genus $m$ Wicks form,
where $g=2m+1$ (see \cite{V2} for details).

\par {\bf Definition 4.} We will say that the word $w^{\prime}$ is obtained from a non-orientable Wicks form $w$ by reduction with respect to the edge $a$.

\section{Proof of the Theorem 1}

We treat the orientable case first.
Since the word $C$ has genus $g$ , it can be obtained from
a Wicks form $U$ of genus $g$ by a non-cancelling substitution
(it means, that all letters of $U$ are replaced by words from a free group such that no cancellations occur). Any Wicks form of genus $g$
can be obtained from the one of maximal length by substituting some letters
by empty words (\cite{V}), that is why
 without loss of generality
we can assume, that $U$ has maximal length $12g-6$ (some letters in the Wicks form can be thought
as replaced by empty words).
By \cite{V}, any Wicks form of maximal length always contains
a negative vertex. If $V$ is a vertex with leaving edges
$a,b,c$, then $$C=u_1a^{-1}bu_2c^{-1}au_3b^{-1}cu_4,$$ where
the vertex $V$ appear as a part of one of the subgraphs on the figure 1
and some of subwords $u_1,u_2,u_3,u_4$ may be empty according to this.
 Now we will
rewrite $C$ using properties of $U$.
$$C=u_1a^{-1}bu_2c^{-1}au_3b^{-1}cu_4=$$
$$u_1a^{-1}bu_2c^{-1}au_3b^{-1}cu_2^{-1}u_2u_4=$$
$$u_1u_3u_3^{-1}a^{-1}bu_2c^{-1}au_3b^{-1}cu_2^{-1}u_2u_4=$$
$$u_1u_3u_3^{-1}a^{-1}bu_2c^{-1}au_3b^{-1}cu_2^{-1}u_3^{-1}u_1^{-1}u_1u_3u_2u_4=$$
$$u_1u_3[u_3^{-1}a^{-1}cu_2^{-1},u_2c^{-1}b]u_3^{-1}u_1^{-1}u_1u_3u_2u_4=$$
$$[u_1u_3u_3^{-1}a^{-1}cu_2^{-1}u_3^{-1}u_1^{-1},u_1u_3u_2c^{-1}bu_3^{-1}u_1^{-1}]u_1u_3u_2u_4=$$
$$[u_1a^{-1}cu_2^{-1}u_3^{-1}u_1^{-1},u_1u_3u_2c^{-1}bu_3^{-1}u_1^{-1}]u_1u_3u_2u_4.$$

Denote
$u_1a^{-1}cu_2^{-1}u_3^{-1}u_1^{-1}$ by $a_1$ and
$u_1u_3u_2c^{-1}bu_3^{-1}u_1^{-1}$ by $b_1$. It is easy to see,
that $|a_1| \leq 2|C|$ and $|b_1| \leq 2|C|$
Now we can see that the word $u_1u_3u_2u_4$ has genus $g-1$ since
it can be obtained by a non-cancelling substitution from a Wicks form $U´$, which
 is obtained from $U$ by reduction with respect to the vertex $V$.

The word $C=u_1a^{-1}bu_2c^{-1}au_3b^{-1}cu_4$ is obtained
from $u_1u_3u_2u_4$ by transformation $\gamma$, if all $|u_1|+|u_4|$, $|u_2|$, $|u_3|$
are not equal 0.
Indeed, the word $u_1u_3u_2u_4$ is corresponding to $w'=w_1xw_3yw_2zw_4$ in the description
of the transformation $\gamma$ as follows: $u_1$ is corresponding to $w_1x_1$, $u_3$ is corresponding
 to $x_2w_3y_1$, $u_2$ is corresponding to $y_2w_2z_1$,  $u_4$ is
corresponding to $z_2w_4$, $a$ is corresponding $a^{-1}$, $b$ is corresponding to $b^{-1}$,
$c$ is corresponding to $c^{-1}$, and $x_1x_2$ is replaced by $x$,$y_1y_2$ is replaced by $y$,
and $z_1z_2$ is replaced by $y$.

$C$ is obtained by transformation $\beta$ if exactly one of $|u_1|+|u_4|$,$|u_2|$,$|u_3|$
is equal to 0, otherwise $C$ is obtained from $u_1u_3u_2u_4$ by transformation
$\alpha$.

In the case of $\beta$ we can assume without loss of generality that $|u_2|=0$,
then $C=u_1abca^{-1}u_3b^{-1}c^{-1}u_4$.
In this case $u_3=u_2^{\prime}d, u_4=d^{-1}u_3^{\prime}$, where $d$ can be empty, and
$C= [u_1abd^{-1}(u_2^{\prime})^{-1},u_1u_2^{\prime}dca^{-1}u_1^{-1}]u_1u_2^{\prime}u_3^{\prime}$.
The word $u_1u_2^{\prime}u_3{\prime}$ corresponds to $w'=w_1xw_2yw_3$  as follows:
 $u_1$ corresponds to $w_1x_1$, $u_2^{\prime}$ corresponds to $x_2w_2y_1$, $u_3^{\prime}$ corresponds to $y_2w_3$,
and $x_1x_2$ is replaced by $x$ and $y_1y_2$ is replaced by $y$.

The easiest case is the transformation $alpha$ and without loss of generality we will assume that
$C=u_1abcdb^{-1}ec^{-1}d^{-1}e^{-1}a^{-1}u_4$ and the negative vertex is the vertex that the edges
$b$ and $d$ are incoming edges, and $c$ is the leaving edge,
$C=[u_1abce^{-1}a^{-1}u_1^{-1},u_1aedb^{-1}a^{-1}u_1^{-1}]u_1u_2$. The word $u_1u_2$ corresponds
to $w'=w_1xw_2x^{-1}w_3$ by $u_1=w_1x_1$, $u_2=x_2w_2x^{-1}w_3$ where $x_1x_2$ is replaced by $x$.

 For details on transformations $\alpha$, $\beta$, $\gamma$ see \cite{V}

Continuing by induction, we get the statement of the Theorem 1 in the orientable case.

Now let $C$ be a non-orientable word of genus $g$.
Then it can be obtained from a non-orientable Wicks form $w$ by a non-cancelling substitution.
We rewrite $C$ using the fact that $w$ has a letter appearing with the same exponent:

$$C=v_1av_2av_3= (v_1av_2v_1^{-1})^2v_1v_2^{-1}v_3.$$ Replacing $v_1av_2v_1^{-1}$ with $x_1$ it is easy to see, that $|x_1| \leq 2|C|$
and $|v_1v_2^{-1}v_3| \leq |C|$. If $v_1v_2^{-1}v_3$ is orientable, then we use the orientable case of the Theorem 1,
which we've just proved. If
$v_1v_2^{-1}v_3$ is non-orientable, then it can be obtained from non-orientable Wicks form $w^{\prime}$ and we proceed as above.

 Let $C$ be a non-orientable word of genus $g$.
It is easy to see, that a product of a square and a commutator can be rewritten as a product of three squares in the following way:
 $x^2[a,b]=(x^2abx^{-1})^2(xb^{-1}a^{-1}x^{-1}a^{-1}x^{-1})^2(xa)^2$.
We have just proved that $C=d_1^2...d_k^2[a_1,b_1]...[a_m,b_m]$,
where $g=k+2m$ and $|a_i|< 2|C|$,$|b_i| < 2|C|$,$|d_i|<2|C|$ for $i=1,...,g$.
Using our rewriting equality, we get that $|a_i| < 6m 2|C| \leq 6g|C| \leq 3|C|^2$.

 Theorem 1 is proved.

\section{Ol'shanskii's result}\label{sec:olshanskii}

The following is proved in \cite{Olshanskii-1989}.

{\bf Theorem 4}
  Let $E$ be a quadratic equation over $F(A)$ in standard form. If
  $g=0,m=2$, or $E$ is not orientable and $g=1,m=1$ then we set
  $N=1$. Otherwise we set $N \leq 3(m-\ol{\chi}(E))$. $E$ has a
  solution if and only if for some $n \leq N$;
  \begin{itemize}
  \item[(i)] there is a set $P = \{p_1,\ldots p_n\}$ of variables and
    a collection of $m$ discs $D_1,\ldots, D_m$ such that,
  \item[(ii)] the boundaries of these discs are circular 1-complexes
    with directed and labelled edges such that each edge has a label in
    $P$ and each $p_j \in P$ occurs exactly twice in the union of
    boundaries;
  \item[(iii)] Let $\ol{\chi}(E)=2-2g$ for orientable surface and $\ol{\chi}(E)=2-g$ for non orientable surface. If we glue the discs together by edges with the same
    label, respecting the edge orientations, then we will have a
    collection $\Sigma_0,\ldots,\Sigma_l$ of closed surfaces and the
    following inequalities: if $E$ is orient able then each $\Sigma_i$
    is orientable and \[\biggl(\sum_{i=0}^{l} \chi(\Sigma_i)\biggr) -
    2l \geq \ol{\chi}(E)\] if $E$ is non-orientable either at least
    one $\Sigma_i$ is non-orientable and \[\biggl(\sum_{i=0}^{l}
    \chi(\Sigma_i)\biggr) - 2l \geq \ol{\chi}(E)\] or, each $\Sigma_i$
    is orientable and
    \[\biggl(\sum_{i=0}^{l} \chi(\Sigma_i)\biggr) - 2l \geq \ol{\chi}(E) +2 \] and
  \item[(iv)] there is a mapping $\ol{\psi}:P \rightarrow (A\cup A^\mo)^*$ such
    that upon substitution, the coefficients $C_1,\ldots,C_{m-1}$ and
    $C$ can be read without cancellations around the boundaries of
    $D_1,\ldots,D_{m-1}$ and $D_m$, respectively; and finally that
  \item[(v)] if $E$ is orientable the discs $D_1,\ldots,D_m$ can be
    oriented so that $w_i$ is read clockwise around $\bdy D_i$ and $d$
    is read clockwise around $\bdy D_m$, moreover all these
    orientations must be compatible with the gluings.
  \end{itemize}

\begin{proof}
  It is shown in Sections 2.4 \cite{Olshanskii-1989} that the
  solvability of a quadratic equation over $F(A)$ coincides with the
  existence of a \emph{diagram $\Delta$ over $F(A)$} on the
  appropriate surface $\Sigma$ with boundary. This diagram may not be
  \emph{simple}, so via surgeries we produce from $\Sigma$ a finite
  collection of surfaces $\Sigma_1,\ldots,\Sigma_l$ with induced
  simple diagrams $\Delta_1,\ldots \Delta_l$ which we can recombine to
  get back $\Sigma$ and $\Delta$. So existence of a diagram $\Delta$
  on $\Sigma$ is equivalent to existence of a collection of simple
  diagrams $\Delta_i$ on surfaces $\Sigma_i$ such that the
  inequalities involving Euler characteristics given in the statement
  of the Theorem are satisfied.

  In Section 2.3 of \cite{Olshanskii-1989} the bounds on $n$ are
  proved. It is also shown in that section that if one can glue discs
  together as described in the statement of the Theorem with the
  condition on the boundaries , then there exist simple diagrams
  $\Delta_i$ on surfaces $\Sigma_i$.
\end{proof}

\section{Quadratic sets of words}

 An {\it orientable quadratic set of words}  is a set of cyclic  words
$ w_1,w_2,\dots ,w_{k}$
 (a cyclic word is the orbit of a linear word under cyclic permutations)
 in some alphabet $a_1^{\pm 1},a_2^{\pm 1},\dots$ of letters
 $a_1,a_2,\dots$ and their inverses $a_1^{-1},a_2^{-1},\dots$ such that
\begin{itemize}
\item[(i)] if $a_i^\epsilon$ appears in $w_i$ (for $\epsilon\in\{\pm 1\}$)
 then $a_i^{-\epsilon}$ appears exactly once in $w_j$,
\item[(ii)] the word $w_i$ contains no cyclic factor (subword of
 cyclically consecutive letters in $w$) of the form $a_l a_l^{-1}$ or
 $a_l^{-1}a_l$ (no cancellation),

\end{itemize}

The {\it genus}  of a quadratic set of words is defined as the sum of genera of the surfaces obtained from
discs $k$ with words $ w_1,w_2,\dots ,w_{k}$ on their boundaries.

{\bf Proof of Theorem 2.}
Let's consider an orientable case first.

Let $W_1$,...,$W_k$ be an orientable quadratic set of words in an alphabet
$\mathcal{A}=\{a_1^{\pm 1},a_2^{\pm 1},\dots\}$
of  genus $g$, and $C_1,\ldots, C_k$ be elements of a free group $F$ such that the system $W_i=C_i,\ i=1,\ldots, k$ has a non-cancelable solution $\phi: F(\mathcal{A})\to F$ in $F$  and $\sum_{i=1}^k|C_i|=s$ in $F$.

We will transform the system $W_i=C_i,\ i=1,\ldots, k$ to one equality.

Suppose there is a letter $a_1$ contained in two equations
 $C_1=W_1=U_1 a_{1}U_2$ and $C_2=W_2=U_3 a_{1}^{-1}U_4$ .
This corresponds to the substitution of $a_1$ from the first equation to the second.
We obtain $1=C_2U_4^{-1}U_1^{-1}C_1U_2^{-1}U_3^{-1}$.  Until there is a letter,
contained in two different equations, we rewrite these two equations as one.
As soon as both appearances of every letter from $\mathcal{A}$ are contained
in one word $L_j=1$,
 we write all $L_j$ next to each other in any order and
 obtain a word $W$ quadratic in $a_1,\ldots,a_n$ containing each $C_i$ only once.
We claim, that some conjugates of $C_i$
in any order
can be presented as a product of at most $g$ commutators  of elements with lengths strictly
less then $2s$, and conjugating elements also have length bounded by $2s$.

Without loss of generality we can assume that the last coefficient
$C_k$ of the standard form appears at the  right end of the word $W$.
 Now let $W=M_3C_{k-1}M_4C_{k}$. Then by conjugating
$C_{k-1}$ by $M_4^{-1}$ we get $W=M_3C_{k-1}^{M_4^{-1}}C_{k}$

The length
of $\phi(W)$ is bounded by $2s$ because the sum of length of all $\phi(a_i)$'s taken twice is not larger than $s$,
and the sum of $C_i's$ is bounded by $s$ as well.
This proves that conjugating element $M_1$ no longer than $2s$.
Continuing by induction, we get that the product of some conjugates of $C_i's$ equals
to a quadratic word  of genus $g$.
(We have genus $g$ since the initial quadratic set of words had genus $g$, this has been proved in the proof of Proposition,
$(a)\rightarrow (b)$).
This word
can be presented as a Wicks form of genus $g$.
Now by using the statement of the Theorem 1, we represent the product of
conjugates as the product of commutators of elements of length not larger than $2s$.

In the non-orientable case there are two cases: when the letter
from $\mathcal{A}$ which we use to bring two equations together,
appears with different exponents or with the same. If a letter
appears in two equations with different exponents, we treat these equations
the same as in the orientable case. But if a letter saying $a_t$
appears in $W_h$ and $W_p$ with same exponent, and $W_h$ and $W_p$ don't contain
a letter with different exponents, namely $W_h=X_1a_tX_2=C_h$
 $W_p=X_3a_tX_4=C_p$, then  the equation obtained  by combining
 $W_h=X_1a_tX_2=C_h$ and
 $W_p=X_3a_tX_4=C_p$ contains $C_p$ and $C_h$ with different exponents.
So, our final word $W$, obtained by bringing all the equations
of the system together, can contain the coefficients with  exponent $-1$,
such that every of these coefficients  $C_h$'s appears as $C_h^{-1}$ in $W$.
Let $a_t$ be a letter with appears with the same exponent in $W$.
By conjugation, we can bring all  $C_h$'s next to $a_t$ such that the length
of  $\phi(W^{\prime})$ is not more then $2s$.
$W^{\prime}=Y_1a_t \tilde{C^{-1}_{h_1}}\tilde{C^{-1}_{h_2}}...\tilde{C^{-1}_{h_r}}Y_2a_tY_3$, where
$\tilde{C^{-1}_{h_1}},\tilde{C_{h_2}},...\tilde{C^{-1}_{h_r}}$
are conjugates of
$C^{-1}_{h_1},C_{h_2},...C^{-1}_{h_r}$.
Without loss of generality we can assume, that $r \leq k/2$, where $k$ is the number
of equations in the system
(otherwise, we take the inverse equation).
Next to the second appearance of $a_t$ we insert $\tilde{C^{-1}_{h_1}}\tilde{C^{-1}_{h_2}}...\tilde{C^{-1}_{h_r}}(\tilde{C^{-1}_{h_1}}\tilde{C^{-1}_{h_2}}...\tilde{C^{-1}_{h_r}})^{-1}$. Substituting $a_t \tilde{C^{-1}_{h_1}}\tilde{C^{-1}_{h_2}}...\tilde{C^{-1}_{h_r}}$ by $a_t^{\prime}$,
we get a new equality, which image under $\phi$ has length not more then $4rs \leq 2ks$
(since $r \leq k/2$) and such that
all coefficients are with the same exponent. Proceeding as in the orientable
case, we get our equality to the standard form, taking into account
the non-orientable case of the Theorem 1.
 First we obtain the product of the conjugates of coefficients in any order
followed by
$d_1^2\ldots d_k^2[a_{1},b_{1}]\ldots[a_m,b_m]$ with $k+2m=g$ and $\phi (a_i), \phi (b_i), \phi (d_i)$ no longer than $4ks$,
and then use the same transformation as in the proof of Theorem 1.
Since $k\leq s$, we get the statement of  Theorem 2.

\section{Hyperbolic groups}

In this section we will prove Theorem 3.

Let $h$ be a word of genus $g$ in $\Gamma=<X|R>$ such that
 $h=[a_1,b_1]...[a_g,b_g]$. The orientable word $U=[A_1,B_1]...[A_g,B_g]$
is a genus $g$ orientable Wicks form.
Let $L=\{A_1,B_1,...,A_g,B_g\}$ and let  $\phi : F(L) \to F(X)$ (where $F(L)$ ($F(X)$) is a free group with basis $L$ ($X$))
be a homomorphism $\phi(A_i)=a_i$, $\phi(B_i)=b_i$ with $i=1,...,g$.
We call this a labelling function for $U$. Let $\mathcal{F}$ be the set
of pairs $(U,\phi)$ where $U$ is a genus $g$ Wicks form and $\phi$ is a labelling
function for $U$ such that $\phi(U)$ is conjugated to $h$ in $\Gamma$.

Consider a pair $(W,\theta)$ in which $|\theta(W)|$ is minimal amongst
all pairs in $\mathcal{F}$. Clearly $\theta(E)$ is minimal for each letter
$E$ in $W$. For convenience we shall take $\hat W$ to be a cyclic permutation
of $W$ such that the last letter of $\hat W$ is labelled by a word
of length more than $12l+M+4$. Consider all letters of $\hat W$ which have
labels greater than $12l+M+4$ in $F(X)$, all  the labels of other letters
are shorter than $12l+M+4$. Corresponding edges in the Cayley graph $G$ of $\Gamma$
will be called long and short edges respectively.

Consider $\theta(\hat W)$ as a path in the Cayley graph $G$. Let $C$
 be  a word in $F(X)$ which represents geodesic for $\theta(\hat W)$ and let $R$ be a minimal word such that $h=_{\Gamma}RCR^{-1}$. We can assume that $2|R|\leq$\{length of all short edges\} (otherwise we can take $\hat W$ to be a cyclic permutation
of $W$ such that the first letter of $\hat W$ is labelled by a word
of length more than $12l+M+4$).  By Lemma 4.3 of \cite{Fulthorp} the terminal vertex of each long edge in $G$
is within $5l+M+3$ of some vertex of $C$, see Figure 2.  ($C$ is corresponding to $F$ in the notation of \cite{Fulthorp}). Let $B$ be a long edge in $\hat W$
which is not
the first long edge in the sequence of letters. Since $\hat W$ is quadratic,
$B$ appears twice, once with exponent 1 and once with exponent -1.
First we shall consider the appearance of $B$ with exponent 1.
In the sequence of letters of $\hat W$, let $A^{\pm 1}$ be the long edge before
$B$ in the sequence such that no long edge appears between $A^{\pm 1}$
and $B$ (note that $A^{\pm 1}$ could be $B^{-1}$).

Let $\iota (p), \tau (p)$ denote the beginning and the end of the path $p$.

\begin{figure}[ht]
\begin{center}
\psfrag{A}{{\scriptsize $\theta(A^{\pm 1})$}}
\psfrag{l}{{\scriptsize $\leq
    5l\hspace{-1pt}+\hspace{-1pt}M\hspace{-1pt}+\hspace{-1pt} 3$}}
\psfrag{B}{{\scriptsize $\theta(B)$}}
\psfrag{u}{{\scriptsize $u$}}
\psfrag{F}{{\scriptsize $C$}}
\psfrag{v}{{\scriptsize $v$}}
\psfrag{all short edges}{all short edges}
\includegraphics[scale=0.4]{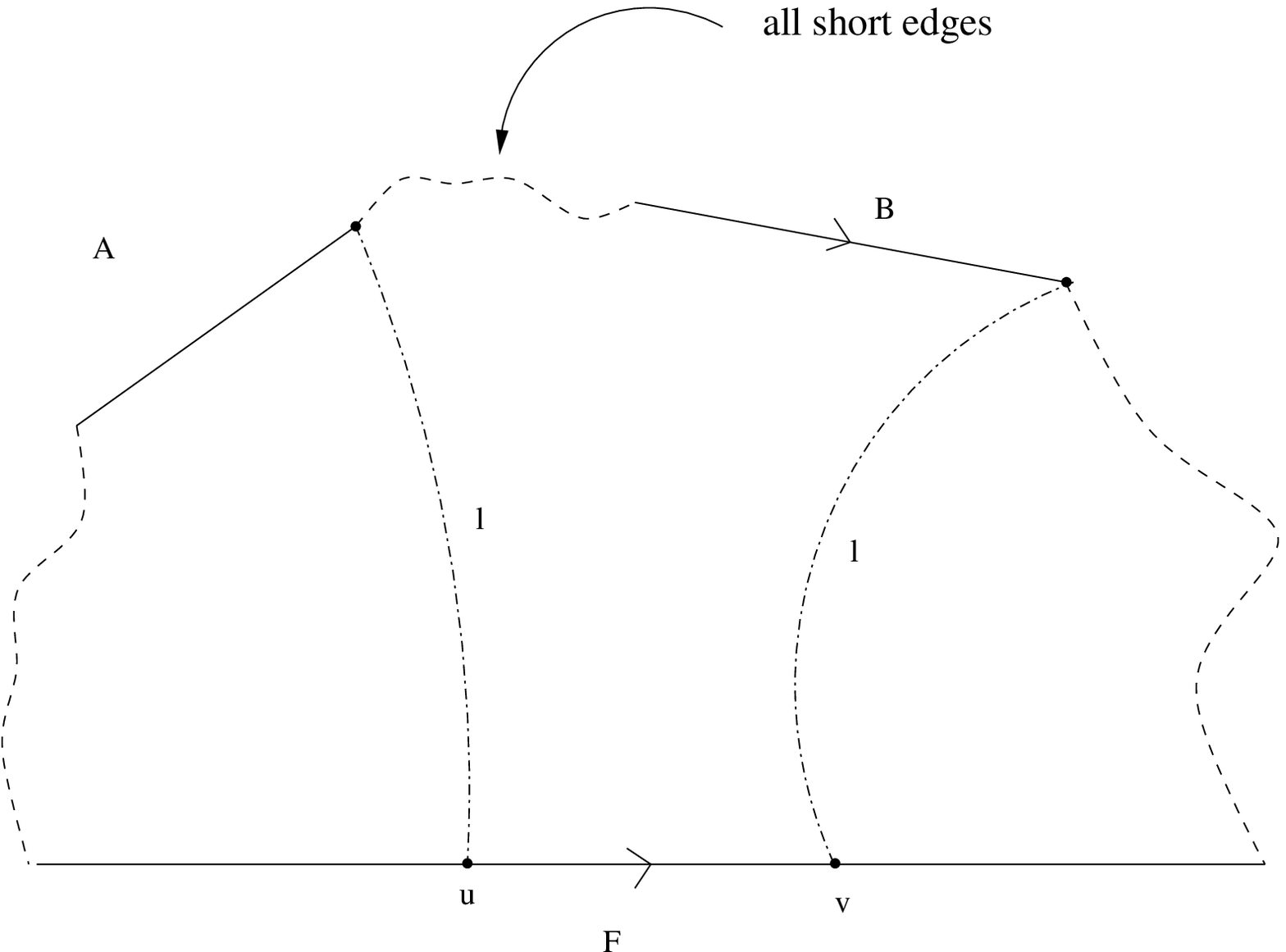}
\\
\vspace{0.38cm}
\it Figure 2.
\end{center}
\end{figure}

{\bf Lemma 1}  [Lemmas 4.3, 4.5 of \cite{Fulthorp}] There exist vertices $u$ and $v$ on $C$, such that
$d(\tau(\theta(A^{\pm1})),u),d(\tau(\theta(B)),v) \leq 5l+M+3$.
We can choose $u$ and $v$ such that
$d(\iota(C),u)<d(\iota(C),v)$.

{\bf Lemma 2} [Follows from the proof of Lemmas 4.3, 4.5 of \cite{Fulthorp}] Let $\hat W=A_1A_2...A_n$ and $A_{i_1},A_{i_2},...,A_{i_k} (i_k=n),$ be long edges.
We can choose vertices $u_1,u_2,...,u_k$ on $C$ such a way that $d(\tau(\theta(A_{i_j})),u_j) \leq 5l+M+3$ and  vertices $u_1,u_2,...,u_k$ appear in the path $C$ in the natural order.

\begin{proof} Lemma 4.3 states the existence of a vertex $u_{j}$ for each long edge $A_{i_j}$ such that $d(\tau(\theta(A_{i_j})),u_{j}) \leq 5l+M+3$. The vertices $u_{j}$ and $u_{j+1}$ chosen in Lemma 4.3 for arbitrary two consecutive long edges, are, actually, chosen such a way that $d(\iota(C),u_j)<d(\iota(C),u_{j+1})$. Indeed,
in the proof of Lemma 4.5 it is shown that the assumption $d(\iota(C),u_j)\geq d(\iota(C),u_{j+1})$ implies a contradiction.\end{proof}
\vspace{2mm}
We now represent accordingly $C$ as $C=D_1...D_k$.

Applying triangle inequalities to long edges, we obtain

\vspace{2mm}
$|\theta A_{j}| \leq 2(5l+M+3)+D_j+$\{length of  short edges between
$A_{j-1}$ and} $A_{j}\}$;

\vspace{2mm}
$\sum|\theta (A_j)| \leq |C| + 2k(5l+M+3)+$\{length of all short edges\};

Therefore

$|\theta(\hat W)| \leq |C| +2k(5l+M+3)+2$\{length of all short edges\}

$\leq
|C|+(12g-6)(12l+M+4).$

Indeed, if $2k_1$ is the number of short edges, then $2k+2k_1\leq 6(2g-1).$
 Notice now that from the triangle inequality,
 $|C|\leq |h|+2|R|\leq|h|+$\{length of all short edges\}. Therefore,

 $|\theta(\hat W)| \leq |h|+(18g-9)(12l+M+4).$

And now we complete the proof using  Theorem 1.

{\bf Acknowledgements}
We are very thankful to the referee for many useful comments and suggestions.


\begin{thebibliography}{99}



\bibitem {Arzh} G.N. Arzhantseva, {\em On quasiconvex subgroups of word hyperbolic groups}, Geometriae dedicata, 87, 191-208, 2001.
\bibitem{BL} J.L.Brenner, R.C.Lyndon, {\em Permutations
 and cubic graphs}, Pacific Journal of Maths, v.104, 285--315

\bibitem{C} M.Culler, {\em Using surfaces to solve equations
 in free groups}, Topology, v.20(2), 1981

\bibitem{CCE} J.A.Comerford, L.P.Comerford and C.C.Edmunds,
 {\em Powers as products of commutators},  Commun. in Algebra, 19(2),
 675--684 (1991)

\bibitem{CE} L.Comerford, C.Edmunds, {\em Products of commutators
 and products of squares in a free group}, Int.J. of Algebra
 and Comput., v.4(3), 469--480, 1994

\bibitem{Comerford-Edmunds-1981}
Leo~P. Comerford, Jr. and Charles~C. Edmunds.
\newblock Quadratic equations over free groups and free products.
\newblock {\em J. Algebra}, 68(2):276--297, 1981.
\bibitem{DGH} V. Diekert, C. Gutierrez, C. Hagenah,
\newblock The existential theory of equations with rational constrains in free groups is PSPACE-complete,
\newblock {\em Information and Computation,} Volume 202, Issue 2 , 1 November 2005, 105-140.

\bibitem{Fulthorp}  S.M. Fulthorp, {\em Squares,Commutators and Genus in Infinite Groups},PhD thesis, School of Mathematics and Statistics, Newcastle University, 2004; {\em Genus $n$ forms over Hyperbolic Groups},
ArXiv:10051513v1[MathGR], 2010.


\bibitem{Grigorchuk-Kurchanov-1989}
R.~I. Grigorchuk and P.~F. Kurchanov.
\newblock On quadratic equations in free groups.
\newblock In {\em Proceedings of the International Conference on Algebra, Part
  1 (Novosibirsk, 1989)}, volume 131 of {\em Contemp. Math.}, pages 159--171,
  Providence, RI, 1992. Amer. Math. Soc.

\bibitem{Grigorchuk-Lysenok-1992}
R.~I. Grigorchuk and I.~G. Lysionok.
\newblock A description of solutions of quadratic equations in hyperbolic
  groups.
\newblock {\em Internat. J. Algebra Comput.}, 2(3):237--274, 1992.
\bibitem{KLMT} O. Kharlampovich, I. Lysenok, A. Myasnikov, N. Touikan, Quadratic
equations over free groups are NP-complete, TOCS (Teor. Comp. Syst.), 10,
2008.
\bibitem{Imp} O. Kharlampovich, A. Myasnikov, \newblock Implicit function theorem over free groups.
\newblock {\em Journal of Algebra},  290 (2005), 1-203.

\bibitem{Malcev-1962}
A.~I. Malcev.
\newblock On the equation {$zxyx\sp{-1}y\sp{-1}z\sp{-1}= aba\sp{-1}b\sp{-1}$}
  in a free group.
\newblock {\em Algebra i Logika Sem.}, 1(5):45--50, 1962.


\bibitem{LS} R.C.Lyndon, P.E.Schupp,{\em
Combinatorial group theory}, Ergebnisse der Mathematik und ihrer Grenzgebiete, Band 89. Springer-Verlag, Berlin-New York, 1977
\bibitem{Mac} J. Macdonald,
\newblock Compressed words and automorphisms in fully residually free groups, \newblock {em International Journal of Algebra and Computation,} 20, no.3 (2010) 343-355.

\bibitem{M} L.Mosher, {\em A User's Guide to the Mapping Class
 Group: Once Punctured Surfaces}, DIMACS Series, Vol.25, 1994.
\bibitem{Olshanskii-1989}
A.~Yu. Olshanskii.
\newblock Diagrams of homomorphisms of surface groups.
\newblock {\em Sibirsk. Mat. Zh.}, 30(6):150--171, 1989.
\bibitem {Olsh} A.~Yu. Olshanskii, \newblock On residualing homomorphisms and $G$-subgroups of hyperbolic groups,
IJAC, 3(4), 1993, 365-409.
\bibitem{V} A.A.Vdovina, {\em Constructing Orientable Wicks Forms
 and Estimation of Their Number}, Communications in Algebra 23 (9),
 3205--3222 (1995).
\bibitem{V2}
A.~Vdovina,
{\em On the number of nonorientable Wicks forms in a free group,}
Proc. R. Soc. Edinb., Sect. A 126, No.1, 113-116 (1996).

\end{thebibliography}
\end{document}